\documentclass[12pt,a4paper, twoside, reqno]{amsart} 
\usepackage{graphicx}
\usepackage[T1]{fontenc} 
\usepackage[latin1]{inputenc} 
\usepackage{bbm} 
\usepackage[english]{babel} 
\usepackage{varioref} 
\usepackage{amsmath,amssymb,mathrsfs,amsthm,amsthm}
\usepackage{setspace} 
\title[Cohomology of Certain Homogeneous Vector Bundles]{On the Cohomology of Certain Homogeneous Vector Bundles of $G/B$ in Characteristic Zero.}
\author{M. Fazeel Anwar}

\begin{document}

\footnotetext[1]{Research supported by COMSATS Institute of Information Technology (CIIT), Islamabad, Pakistan.}

\begin{abstract}
\noindent In his famous paper \cite{Demazure}, Demazure introduced certain indecomposable modules and used them to give a short proof of Bott's theorem. In this paper we consider a generalization of these modules and give their cohomology. 
\end{abstract}

\maketitle 

\noindent Let $G$ be a reductive connected algebraic group over an algebraically closed field $k$ of characteristic zero and $B$ be a Borel subgroup. In his famous paper \cite{Demazure}, Demazure introduced the indecomposable modules $V_{\lambda,\alpha}$ with weights $\lambda,\lambda-\alpha,...,s_{\alpha}(\lambda)$, where $\alpha$ is a simple root and  $s_{\alpha}$ is the corresponding reflection. He  used these modules to give a short proof of Bott's theorem. In this paper, we consider a generalization of the module $V_{\lambda,\alpha}$ denoted by $M_{\alpha,r}(\lambda)$. This is the unique (up to isomorphism) indecomposable $B-$module with weights $\lambda,\lambda-\alpha,...,\lambda-r\alpha$. We determine the $i-$th cohomology of $M_{\alpha,r}(\lambda)$ for all $i$. This in particular gives all cohomology of the modules $V_{\lambda,\alpha}$. 

\bigskip

\noindent Let $T \subset B$ be a maximal torus of $G$. Let ${\rm mod}(G)$ be the category of finite dimensional rational $G-$modules. Define $X(T)$ to be the group of multiplicative characters of $T$. The Weyl group $W$ acts on $T$ in the usual way. Let $\Phi$ be the system of roots, $\Phi^{+}$ the set of positive roots for which $B$ is the negative Borel subgroup and $S$ the set of simple roots. For $\alpha \in \Phi$ its coroot $\alpha^{v}$ is given by $\frac{2\alpha}{(\alpha,\alpha)}$, where $(\, , )$ is a positive definite, $W-$invariant, symmetric, non-singular, bilinear form on $\mathbb{R}\otimes_{\mathbb{Z}} X(T)$. For $\alpha \in S$, the parabolic subgroup of $G$ containing $B$ is denoted by $P_{\alpha}$ and has $\alpha$ as its only positive root.   

\bigskip

\noindent We refer the reader to \cite{Jantzen},\cite{Fazeel} and \cite{Donkin} for terminology and results not explained here.

\bigskip

\noindent We first consider the module $M_{\alpha,1}(\lambda)$. This module is denoted by $N_{\alpha}(\lambda)$ in \cite{Donkin2007} and plays a crucial role in finding a recursive description for the characters of cohomology of line bundles on three dimensional flag variety(in characteristic p).  

\bigskip

\noindent {\bf Theorem 1:} Let $\alpha$ be a simple root and $\lambda \in X(T)$ then

$$H^{i}(M_{\alpha,1}(\lambda))=\begin{cases}
H^{i}(\lambda) \oplus H^{i}( \lambda-\alpha), &\langle \lambda , \alpha^{v} \rangle \neq 0\\
0, &\langle \lambda , \alpha^{v} \rangle = 0.
\end{cases}$$

\bigskip

\noindent {\bf Proof:} We will give the proof in separate cases.

\noindent {\bf a:} Let $\langle \lambda , \alpha^{v} \rangle \leq -1$. On the second page of the spectral sequence we have $$H^{i}(M_{\alpha,1}(\lambda))= R^{i-1}{\rm Ind}_{P_{\alpha}}^{G} R^{1}{\rm Ind}_{B}^{P_{\alpha}}(M_{\alpha,1}(\lambda)).$$ 
Also since $\langle \lambda , \alpha^{v} \rangle \leq -1$ so ${\rm Ind}_{B}^{P_{\alpha}}\lambda=0$. Moreover $P_{\alpha}/B$ is one dimensional so $R^{i}{\rm Ind}_{B}^{P_{\alpha}}\lambda=0$ for all $i \geq 2$. Hence from the short exact sequence 
$$0 \rightarrow \lambda-\alpha \rightarrow M_{\alpha,1}(\lambda) \rightarrow \lambda \rightarrow 0$$ 
we get $$0 \rightarrow R^{1}{\rm Ind}_{B}^{P_{\alpha}}(\lambda-\alpha) \rightarrow R^{1}{\rm Ind}_{B}^{P_{\alpha}}(M_{\alpha,1}(\lambda)) \rightarrow R^{1}{\rm Ind}_{B}^{P_{\alpha}}(\lambda) \rightarrow 0.$$

\noindent Since all modules for $P_{\alpha}/R_{u}(P_{\alpha})$ are completely reducible so we get $R^{1}{\rm Ind}_{B}^{P_{\alpha}}(M_{\alpha,1}(\lambda)) \simeq R^{1}{\rm Ind}_{B}^{P_{\alpha}}(\lambda-\alpha) \oplus R^{1}{\rm Ind}_{B}^{P_{\alpha}}(\lambda)$. Therefore $$H^{i}(M_{\alpha,1}(\lambda))= R^{i-1}{\rm Ind}_{P_{\alpha}}^{G}(R^{1}{\rm Ind}_{B}^{P_{\alpha}}(\lambda-\alpha) \oplus R^{1}{\rm Ind}_{B}^{P_{\alpha}}(\lambda))$$ 
and we get the result. 

\bigskip

\noindent {\bf b:} For $\langle \lambda,\alpha^v \rangle = 1$ we get $R^{j}{\rm Ind}_{B}^{P_{\alpha}}(\lambda-\rho)$ is zero for all $j \neq 0$. Therefore $\nabla_{\alpha}(\rho) \otimes \nabla_{\alpha}(\lambda-\rho) =\nabla_{\alpha}(\lambda)$. Hence $H^{i}(M_{\alpha,1}(\lambda))= H^{i}(\lambda)$.

\bigskip

\noindent {\bf c:} Suppose $\langle \lambda,\alpha^v \rangle \geq 2$. On the second page of the spectral sequence we have $$R^{i}{\rm Ind}_{P_{\alpha}}^{G}R^{j}{\rm Ind}_{B}^{P_{\alpha}}(\nabla_{\alpha}(\rho)\otimes(\lambda-\rho))=R^{i}{\rm Ind}_{P_{\alpha}}^{G}(\nabla_{\alpha}(\rho)\otimes R^{j}{\rm Ind}_{B}^{P_{\alpha}}(\lambda-\rho)).$$ For $\langle \lambda,\alpha^v \rangle \geq 2$ we have $R^{j}{\rm Ind}_{B}^{P_{\alpha}}(\lambda-\rho)$ is zero for all $j \neq 0$. Therefore $H^{i}(M_{\alpha,1}(\lambda))= R^{i}{\rm Ind}_{P_{\alpha}}^{G}(\nabla_{\alpha}(\rho) \otimes \nabla_{\alpha}(\lambda-\rho))$. Since the weights of $\nabla_{\alpha}(\rho)$ are $\rho$ and $\rho-\alpha$, by a special case of Clebsch Gordan formula we get    $\nabla_{\alpha}(\rho) \otimes \nabla_{\alpha}(\lambda-\rho) \simeq \nabla_{\alpha}(\lambda) \oplus \nabla_{\alpha}(\lambda-\alpha)$ (we are working in characteristic zero) and hence   $$H^{i}(M_{\alpha,1}(\lambda))= H^{i}(\lambda) \oplus H^{i}(\lambda-\alpha).$$

\bigskip

\noindent {\bf d:} Finally we consider the case $\langle \lambda,\alpha^v \rangle = 0$. We have $\langle \lambda-\rho,\alpha^v \rangle = -1$ and hence $R^{j}{\rm Ind}_{B}^{P_{\alpha}}(\lambda-\rho)=0$. So we get $H^{i}(M_{\alpha,1}(\lambda))=0$.

\noindent This completes the proof.

\bigskip

\noindent {\bf Theorem 2:} Let $r \geq 0$, $\lambda \in X(T)$ and $s=\langle \lambda-r\rho,\alpha^{v}\rangle$ then

$$H^{i}(M_{\alpha,r}(\lambda))=\begin{cases}
\bigoplus_{t=0}^{r} H^{i}(\lambda-t\alpha), &\langle \lambda , \alpha^{v} \rangle \leq -1 \\ 
\bigoplus_{t=0}^{r} H^{i}(\lambda-t\alpha), &\langle \lambda , \alpha^{v} \rangle > r \\& \text{and} \, r \leq s\\
\bigoplus_{t=0}^{s} H^{i}(\lambda-t\alpha), &\langle \lambda , \alpha^{v} \rangle > r \\& \text{and} \, r>s\\
\bigoplus_{t=0}^{r} H^{i}(t\rho-(t-2-m)\alpha), &0 \leq \langle \lambda , \alpha^{v} \rangle <r-1 \\
H^{i}(\lambda), &\langle \lambda , \alpha^{v} \rangle = r\\ 
0, &\langle \lambda , \alpha^{v} \rangle = r-1.
\end{cases}$$

\bigskip

\noindent {\bf Proof:} We will use mathematical induction on $r$ to prove the result. The result is true for $r=1$ by theorem 1. Suppose the result is true for $r-1$ then 

$$H^{i}(M_{\alpha,r-1}(\lambda))=\begin{cases}
\bigoplus_{t=0}^{r-1} H^{i}(\lambda-t\alpha), &\langle \lambda , \alpha^{v} \rangle \leq -1 \\ 
\bigoplus_{t=0}^{r-1} H^{i}(\lambda-t\alpha), &\langle \lambda , \alpha^{v} \rangle > r-1 \\& \text{and} \, r-1 \leq s\\
\bigoplus_{t=0}^{s} H^{i}(\lambda-t\alpha), &\langle \lambda , \alpha^{v} \rangle > r-1 \\& \text{and} \, r-1>s\\
\bigoplus_{t=0}^{r-1} H^{i}(t\rho-(t-2-m)\alpha), &0 \leq \langle \lambda , \alpha^{v} \rangle <r-2 \\
H^{i}(\lambda), &\langle \lambda , \alpha^{v} \rangle = r-1\\ 
0, &\langle \lambda , \alpha^{v} \rangle = r-2.
\end{cases}$$
\noindent Now for $r$ we give the result in cases as in theorem 1.

\noindent {\bf a:} Let $\langle \lambda , \alpha^{v} \rangle \leq -1$. On the second page of the spectral sequence we have $$H^{i}(M_{\alpha,r}(\lambda))= R^{i-1}{\rm Ind}_{P_{\alpha}}^{G} R^{1}{\rm Ind}_{B}^{P_{\alpha}}(M_{\alpha,r}(\lambda)).$$ 
Moreover we have the short exact sequence $$0 \rightarrow M_{\alpha,r-1}(\lambda-\alpha) \rightarrow M_{\alpha,r}(\lambda) \rightarrow \lambda \rightarrow 0.$$ Also since $\langle \lambda , \alpha^{v} \rangle \leq -1$ so ${\rm Ind}_{B}^{P_{\alpha}}(\lambda)=0$. Moreover $P_{\alpha}/B$ is one dimensional so $R^{i}{\rm Ind}_{B}^{P_{\alpha}}(\lambda)=0$ for all $i \geq 2$. Using the above short exact sequence we get $$0 \rightarrow R^{1}{\rm Ind}_{B}^{P_{\alpha}}(M_{\alpha,r-1}(\lambda-\alpha)) \rightarrow R^{1}{\rm Ind}_{B}^{P_{\alpha}}(M_{\alpha,r}(\lambda)) \rightarrow R^{1}{\rm Ind}_{B}^{P_{\alpha}}(\lambda) \rightarrow 0.$$

\noindent Since all modules for $P_{\alpha}/R_{u}(P_{\alpha})$ are completely reducible so we get $$R^{1}{\rm Ind}_{B}^{P_{\alpha}}(M_{\alpha,r}(\lambda)) \simeq R^{1}{\rm Ind}_{B}^{P_{\alpha}}(M_{\alpha,r-1}(\lambda-\alpha)) \oplus R^{1}{\rm Ind}_{B}^{P_{\alpha}}(\lambda).$$ 
\noindent Therefore  $$H^{i}(M_{\alpha,r}(\lambda))= R^{i-1}{\rm Ind}_{P_{\alpha}}^{G}(R^{1}{\rm Ind}_{B}^{P_{\alpha}}(M_{\alpha,r-1}(\lambda-\alpha)) \oplus R^{1}{\rm Ind}_{B}^{P_{\alpha}}(\lambda)).$$ Since we are working in characteristic zero we can get $$H^{i}(M_{\alpha,r}(\lambda))=H^{i}(M_{\alpha,r-1}(\lambda))\oplus H^{i}(\lambda).$$ 
Now use the inductive hypothesis to get the result. 

\bigskip

\noindent {\bf b:} Suppose $\langle \lambda,\alpha^v \rangle \geq r$. On the second page of the spectral sequence we have 
$$R^{i}{\rm Ind}_{P_{\alpha}}^{G}R^{j}{\rm Ind}_{B}^{P_{\alpha}}(\nabla_{\alpha}(r\rho)\otimes(\lambda-r\rho))=R^{i}{\rm Ind}_{P_{\alpha}}^{G}(\nabla_{\alpha}(r\rho)\otimes R^{j}{\rm Ind}_{B}^{P_{\alpha}}(\lambda-r\rho)).$$ 
For $\langle \lambda,\alpha^v \rangle \geq r$ we have $R^{j}{\rm Ind}_{B}^{P_{\alpha}}(\lambda-r\rho)$ is zero for all $j \neq 0$. Therefore 
$$R^{i}{\rm Ind}_{P_{\alpha}}^{G}R^{j}{\rm Ind}_{B}^{P_{\alpha}}(M_{\alpha,r}(\lambda))= R^{i}{\rm Ind}_{P_{\alpha}}^{G}(\nabla_{\alpha}(r\rho) \otimes \nabla_{\alpha}(\lambda-r\rho)).$$

\noindent Now we have two cases here. Firstly let $r \leq s$ then we will get $H^{i}(M_{\alpha,r}(\lambda))= \bigoplus_{t=0}^{r} H^{i}(\lambda-t\alpha)$. Now if $r>s$ then we have $$H^{i}(M_{\alpha,r}(\lambda))= \bigoplus_{t=0}^{s} H^{i}(\lambda-t\alpha).$$

\bigskip

\noindent {\bf c:} For $0 \leq \langle \lambda,\alpha^v \rangle < r-1$ we get $R^{j}{\rm Ind}_{B}^{P_{\alpha}}(\lambda-r\rho)$ is zero for all $j \neq 1$. We get $$R^{i}{\rm Ind}_{P_{\alpha}}^{G}R^{j}{\rm Ind}_{B}^{P_{\alpha}}(\nabla_{\alpha}(r\rho)\otimes(\lambda-r\rho))$$ $$=R^{i}{\rm Ind}_{P_{\alpha}}^{G}(\nabla_{\alpha}(r\rho)\otimes R^{1}{\rm Ind}_{B}^{P_{\alpha}}(\lambda-r\rho)).$$ 

\noindent Using Serre duality we get $(R^{1}{\rm Ind}_{B}^{P_{\alpha}}(\lambda-r\rho))^{*}={\rm Ind}_{B}^{P_{\alpha}}(-\lambda+r\rho-\alpha)$. Therefore $H^{i}(M_{\alpha,r}(\lambda))=R^{i}{\rm Ind}_{P_{\alpha}}^{G}(\nabla_{\alpha}(r\rho)\otimes \nabla_{\alpha}(-\lambda+r\rho-\alpha)^{*})$.

\noindent Now let $\langle \lambda,\alpha^v \rangle=m$ then we get $$H^{i}(M_{\alpha,r}(\lambda))=R^{i}{\rm Ind}_{P_{\alpha}}^{G}(\nabla_{\alpha}(r\rho)\otimes \nabla_{\alpha}(\lambda-r\rho+(r-1-m)\alpha)).$$ 
Apply the Clebsch Gordan formula again to get the result. 

\bigskip

\noindent {\bf d:} For $\langle \lambda,\alpha^v \rangle = r$ we get $R^{j}{\rm Ind}_{B}^{P_{\alpha}}(\lambda-r\rho)$ is zero for all $j \neq 0$. Therefore $\nabla_{\alpha}(r\rho) \otimes \nabla_{\alpha}(\lambda-r\rho) =\nabla_{\alpha}(\lambda)$. Hence $H^{i}(M_{\alpha,r}(\lambda))= H^{i}(\lambda)$.

\bigskip

\noindent {\bf e:} Finally we consider the case $\langle \lambda,\alpha^v \rangle = r-1$ so $\langle \lambda-r\rho,\alpha^v \rangle = -1$ and hence $R^{j}{\rm Ind}_{B}^{P_{\alpha}}(\lambda-\rho)=0$. So we get $H^{i}(M_{\alpha,r}(\lambda))=0$.

\noindent This completes the proof.

\bigskip

\noindent \textit{Acknowledgments}. I am very grateful to Stephen Donkin for bringing this problem to my attention and for his valuable remarks.

\noindent Department of Mathematics, University of York, Heslington, York, YO10
5DD, United Kingdom.\\ 
E-mail: mfa501@york.ac.uk

\end{document}